\documentclass[11pt]{elsarticle}
\usepackage{xcolor}
\usepackage{graphicx} % Required for inserting images
\usepackage{fullpage}
\usepackage{amsmath,amssymb}
\usepackage{hyperref}

\usepackage{bm}
\usepackage{algorithm}%
\usepackage{algorithmicx}%
\usepackage{algpseudocode}%
\usepackage{listings}%

\def\vs{\vspace{0.2cm}}

\def\I{\mathcal{I}}

\def\J{\mathcal{J}}

\def\DEIM{\texttt{DEIM}}
\def\diag{\texttt{diag}}

\begin{document}

\title{Interpolatory dynamical low-rank approximation for the 3+3d Boltzmann--BGK equation}
%\author{Alec Dektor, Lukas Einkemmer}
\author[lbnl]{Alec Dektor}
\author[uibk]{Lukas Einkemmer\corref{cor1}} \ead{lukas.einkemmer@uibk.ac.at}

\address[lbnl]{Applied Mathematics and Computational Research Division, Lawrence Berkeley National Laboratory, Berkeley, CA, USA}
\address[uibk]{Department of Mathematics, Universit\"at Innsbruck, Austria}
\cortext[cor1]{Corresponding author}

\begin{abstract}
We introduce two novel interpolatory dynamical low-rank (DLR) approximation methods for the efficient time integration of the Boltzmann–BGK equation. Both methods overcome limitations of classic DLR schemes based on orthogonal projections for nonlinear equations. In particular, we demonstrate that the proposed methods can efficiently compute solutions to the full Boltzmann–BGK equation without restricting to e.g.~weakly compressible or isothermal flow. The first method we propose directly applies the recently developed interpolatory projector-splitting scheme on low-rank matrix manifolds. The second method is a variant of the rank-adaptive basis update and Galerkin scheme, where the Galerkin step is replaced by a collocation step, resulting in a new scheme we call basis update and collocate (BUC). Numerical experiments in both fluid and kinetic regimes demonstrate the performance of the proposed methods. In particular we demonstrate that the methods can be used to efficiently compute low-rank solutions in the six-dimensional (three spatial and three velocity dimensions) setting on a standard laptop. 
\end{abstract}

\maketitle

\section{Introduction}

The Boltzmann equation is a fundamental model in gas dynamics, providing a statistical description of the positions and velocities of gas molecules. It is particularly effective in regimes where traditional fluid models break down, such as in the study of rarefied gases encountered during atmospheric re-entry. 
Numerical simulations of the full Boltzmann equation are extremely expensive due to the curse of dimensionality (the probability density function that is advanced in time is up to six-dimensional) and the form of the collision operator that appears in the Boltzmann equation (which requires evaluating up to six-dimensional integrals for each point in space and velocity). While the latter challenge can often be mitigated by employing simplified collision models, such as the BGK (Bhatnagar–Gross–Krook) or Fokker–Planck models, the high dimensionality is an inherent feature of kinetic models. When the system is in thermodynamic equilibrium, the Boltzmann equation reduces to the three-dimensional Navier–Stokes equations, which provide accurate approximations in many scenarios. However, for many important applications — such as rarefied gas flow or plasma dynamics — this simplification is not valid.

Recently, dynamical low-rank (DLR) approximations \cite{Koch2007,Lubich2014} have been used successfully for kinetic equations ranging from radiative transfer/transport \cite{Peng2020,Ding2021,Peng2021,Einkemmer2021,Kusch2021a} and Boltzmann equations \cite{Einkemmer2019,Einkemmer2021d,Hu2022} to the Vlasov equation \cite{Kormann2015,Einkemmer2020,Einkemmer2018,Einkemmer2021a,Einkemmer2022,Guo2024,AllmannRahn2022,Guo2024b,Coughlin2024}. Low-rank approximation alleviates the curse of dimensionality by representing the six-dimensional solution as a linear combination of products of at most three-dimensional basis functions. The DLR approach derives equations that govern the time evolution of both the basis functions and their coefficients. For many applications, a small set of basis functions is sufficient to accurately represent the solution at each time step. This has been demonstrated in scenarios where the solution remains close to thermodynamic equilibrium \cite{Einkemmer2021,Einkemmer2024a,Patwardhan2024}, but also in the fully kinetic regime \cite{Einkemmer2018,Peng2020,Einkemmer2020,AllmannRahn2022,Einkemmer2023}. 

For the Boltzmann equation with BGK collision operator (in the following called the Boltzmann--BGK equation for brevity) a number of methods have been proposed in the literature. The main challenge in this case is that the collision operator directly includes the Maxwellian (the equilibrium distribution) which can not be represented easily in a low-rank format. Two methods have been proposed in the literature to address this.
\begin{enumerate}
    \item In the weakly compressible regime a representation can be obtained by expanding the Maxwellian in terms of the velocity divided by the speed of sound (a small quantity). Such an approach has been considered in \cite{Einkemmer2019}.
    \item In \cite{Einkemmer2021d} it was proposed to subject only the multiplicative deviation from equilibrium to a low-rank approximation. This representation also has the advantage that a smaller number of basis functions can be used (see section \ref{sec:model} for more details). A fast convolution algorithm can then be used to compute the required quantities in the algorithm in an efficient way (even though the Maxwellian is not in low-rank form). However, this only works in the isothermal case (i.e.~no or small variations in the temperature of the fluid).
\end{enumerate}
However, currently, it is not possible to treat the general case (which allows for compression as well as temperature variations). The fundamental issue is that the BGK model, because of its dependence on the Maxwellian, is a nonlinear operator that can not be written in low-rank form. The classical DLR schemes based on orthogonal projections applied to this problem are memory efficient (i.e.~the memory cost scales as $n^3$ instead of $n^6$, where $n$ is the number of grid points in each dimension; here assumed to be equal) but not computationally efficient (i.e.~computational cost still scales as $n^6$ in the worst case). 

More recently, DLR schemes based on interpolatory projections have been introduced for integrating low-rank solutions to nonlinear PDEs \cite{Donello2023,Naderi2023,Ghahremani2024a,Ghahremani2024b,Dektor2024}. The advantage of these schemes is that the resulting evolution equations for the DLR basis functions and coefficients only require the evaluation of the vector field defining the PDE at a small set of interpolation points. This approach can therefore efficiently integrate any PDE defined by a vector field that can be efficiently evaluated pointwise. 

In this paper, we propose two solvers for the Boltzmann-BGK equation \eqref{eq:boltzmann-bgk} based on the interpolatory DLR approach. The first solver is obtained by applying the interpolatory projector-splitting scheme recently introduced in \cite{Dektor2024} to the Boltzmann-BGK equation. The second solver is a novel rank-adaptive interpolatory DLR scheme, which is a modification of the basis update and Galerkin (BUG) method introduced in \cite{Ceruti2022} and improved in \cite{Ceruti2023}. In this approach, we modify the basis update steps and replace the Galerkin step with a collocation step, resulting in a scheme we term basis update and collocate (BUC). Both solvers rely on efficient pointwise evaluations of the vector field defining the Boltzmann-BGK equation. We demonstrate that the vector field defining the Boltzmann-BGK equation can indeed be efficiently evaluated pointwise, leading to solvers with computational complexity scaling as $n^3$ instead of $n^6$. 

The remainder of the paper is structured as follows. In Section~\ref{sec:model} we recall the Boltzman-BGK equation and the DLR methodology. We recapitulate the classical DLR approach based on orthogonal projections and discuss its computational challenges for the Boltzmann-BGK equation. In Section~\ref{sec:interp_DLR} we present the interpolatory projector-splitting DLR solver and the basis update and collocate (BUC) solver for the Boltzmann-BGK equation and show that they can efficiently integrate low-rank solutions. In Section~\ref{sec:numerics} we provide several numerical examples demonstrating the efficacy of the proposed solvers in both the fluid and kinetic regimes. Our main findings are summarized in Section~\ref{sec:conclusions}. 

\section{Boltzmann--BGK equation \label{sec:model}}

The Boltzmann--BGK equation is stated as follows
\begin{equation} \label{eq:boltzmann-bgk}
\partial_{t}f+v \cdot \nabla_{x}f=\frac{1}{\epsilon}(M(f)-f),
\end{equation}
where $f(t,x,v)$ is the sought-after particle density function (i.e.~the unknown). Note that discretizing this six-dimensional function is extremely expensive. Choosing $200$ grid points in each direction requires $512$ Terabytes of memory just to store a single copy of $f$. The physical domain is $\Omega = \Omega_x \times \Omega_v$, where $\Omega_x$ is usually bounded and $\Omega_v = \mathbb{R}^3$. In numerical simulations, $\Omega_v$ is often truncated to a box. The strength of the collisionality is determined by $\epsilon$; physically, this is the Knudsen number, which is given by the ratio between the molecular mean free path and the length scale of the fluid.

The BGK collision operator on the right-hand side is a nonlinear function due to nonlinear dependence of the Maxwellian $M(f)$ on $f$. The Maxwellian is given by
\begin{equation} \label{eq:maxwellian}
  M(f)  =\frac{\rho(f)}{(2\pi T(f))^{d_v/2}}\exp\left(-\frac{\| v-u(f) \|^{2}}{2T(f)}\right)
\end{equation}
and is completely determined by the first three moments of the distribution function
\begin{equation}
\rho(f)  =\int_{\Omega_v}f\,dv, \qquad
u(f)  =\frac{1}{\rho(f)}\int_{\Omega_v}vf\,dv,  \qquad
T(f)  =\frac{1}{{\color{black} d_v} \rho(f)}\int_{\Omega_v} \| v-u(f) \|^{2}f\,dv.
\end{equation}
Usually, the Boltzmann equation is used in non-dimensionalized form (as we have done here) such that the speed of sound is $1$. Then $u$ is the velocity in units of the speed of sound (i.e.~the local Mach number of the flow). The dynamical low-rank work in \cite{Einkemmer2019} considers the case where $u \ll 1$. Then the Maxwellian can be expanded in $u$ and the result turns out to be a low-rank approximation of rank $10$ up to $\mathcal{O}(u^3)$. Such expansion have been used extensively e.g.~in lattice Boltzmann methods \cite{Chen1998}. On the other hand, the work in \cite{Einkemmer2021d} considers the case where $T=\text{const}$, but $u$ can be large. We emphasize that in the present paper, no assumptions on the moments are made. 

A particularly interesting limit is $\epsilon \to 0$ (i.e.~the limit of large collisionality). In this case $f$ is a Maxwellian with the moments determined by the Euler equations of fluid dynamics. If terms up to $\mathcal{O}(\epsilon)$ are kept (this can be done e.g.~by a Chapman--Enskog expansion) we still have a Maxwellian distribution function (up to $\mathcal{O}(\epsilon))$. However, the moments are now described by the Navier--Stokes equations with a viscosity that is proportional to $\mathcal{O}(\epsilon)$. For more details we refer the reader to \cite{Bardos1991,Bardos1993}.

\subsection{Dynamical low-rank approximation \label{sec:DLR}}

Before describing the interpolatory DLR approach we recapitulate the classical DLR approach based on orthogonal projections and discuss the associated computational challenges for the Boltzmann-BGK equation \eqref{eq:boltzmann-bgk}. 
At each time $t$ we seek an approximate distribution function {\color{black} $f_r(t)$ in the set 
\begin{equation} \label{eq:low_rank_manifold} 
    \mathcal{M}_r = \{g(x,v) \in 
    L^2(\Omega_x \times \Omega_v) \ : \ g(x,v) = \sum_{i,j=1}^r X_i(x)S_{ij}V_{j}(v), \ S \in GL_r(\mathbb{R}) \}
\end{equation} 
of rank-$r$ functions. Thus our approximate distribution function takes 
}
the form 
\begin{equation} \label{eq:DLR}
f_{\color{black}_r}(t,x,v) = \sum_{i,j=1}^r X_i(t,x) S_{ij}(t) V_j(t,v). 
\end{equation} 
It is convenient to collect the {\color{black} time-dependent} component functions into {\color{black} $r$-dimensional} row vectors $X(t,x) = \left(X_1(t,x),\ldots,X_r(t,x) \right)$ and $V(t,v) = \left(V_1(t,v),\ldots,V_r(t,v) \right)$ so that 
\begin{equation} \label{eq:DLR_matvec}
    f_{\color{black}_r}(t,x,v) = X(t,x) S(t) V(t,v)^{\top}. 
\end{equation}
{\color{black} It is well-known that the set \eqref{eq:low_rank_manifold} is a smooth submanifold of $L^2(\Omega_x \times \Omega_v)$, see, e.g., \cite[Section 3]{dektor2021dynamic}. Thus, for each $f_r(t)$ one has a tangent space containing the infinitesimal perturbations along which $f_r(t)$ remains on the rank-$r$ manifold $\mathcal{M}_r$. The approximate solution $f_r(t)$ can be integrated on $\mathcal{M}_r$ by selecting $\partial_t f_r(t)$ in the tangent space at each time $t$. The classic DLR method selects $\partial_t f_r(t)$ by solving a variational problem at each time $t$ that minimizes the $L^2$ distance between the tangent space and the time-derivative defined by the Boltzmann-BGK equation \eqref{eq:boltzmann-bgk} 
}
\begin{equation} \label{eq:h_def}
h(t,x,v) = -v\partial_{x}f_{\color{black} r}+\frac{1}{\epsilon}(M(f_{\color{black} r})-f_{\color{black} r}). 
\end{equation}
{\color{black} Since $\mathcal{M}_r$ is an embedded submanifold of $L^2(\Omega_x \times \Omega_v)$ it follows that its tangent spaces are vector subspaces of $L^2(\Omega_x \times \Omega_v)$ and the variational problem at time $t$ is solved by projecting \eqref{eq:h_def}} orthogonally onto the tangent space. 
{\color{black} 
To obtain such orthogonal projection it is convenient to select the component functions} so that $X$ and $V$ are orthonoromal\footnote{\color{black} Orthonormal component functions can always be obtained, e.g., using the Gram-Schmidt procedure.}, i.e., 
\begin{equation} \label{eq:orth_bases}
\left\langle X^{\top} X \right\rangle_x = \left\langle V^{\top} V \right\rangle_v = I_{r \times r}, 
\end{equation}
where $\langle g(x) \rangle_x = \int_{\Omega_x} g(x) dx$ and $\langle g(v) \rangle_v = \int_{\Omega_v} g(v) dv$. {\color{black} When the component functions are orthonormal \eqref{eq:orth_bases}, the orthogonal projection of \eqref{eq:h_def} onto the tangent space, and hence the classical DLR evolution equation for $f_r$, is } 
\begin{equation} \label{eq:orth_proj}
\partial_t f_{\color{black} r}(t,x,v) = \left\langle h V \right\rangle_v V^{\top} -  X \left\langle X^{\top} h V \right\rangle_{x,v} V^{\top} + X \left\langle X^{\top} h \right\rangle_{x}. 
\end{equation}
{\color{black} It is well-known that the curvature of the manifold $\mathcal{M}_r$ is inversely proportional to the smallest singular values of the coefficient matrix $S$, which makes the evolution equation \eqref{eq:orth_proj} extremely stiff when $S$ contains singular values near zero. It was demonstrated in \cite{Lubich2014} that operator splitting methods are robust to the curvature of $\mathcal{M}_r$ and have since become widely used in DLR schemes. The simplest of these operator splitting schemes is the first-order Lie-Trotter splitting. 
Applied to \eqref{eq:orth_proj},} Lie-Trotter yields a numerical scheme with the following three substeps (for more details we refer to {\color{black} \cite{Lubich2014}} and \cite{Einkemmer2018,Einkemmer2021d} which use a very similar notation): We start from a rank-$r$ representation of the distribution function \eqref{eq:DLR} with orthonormal bases \eqref{eq:orth_bases}. First, let $K(t,x) = X(t,x) S(t)${\color{black} , a length-$r$ row vector of non-orthogonal component functions, } so that $f_{\color{black} r} =  K V^{\top}$, and update $K$ by integrating 
\begin{equation} \label{eq:K_step_orth}
\partial_t K = \left\langle h V \right\rangle_v. 
\end{equation}
Then perform an orthonormalization of the updated $K$ to obtain {\color{black} an updated set of $r$ orthonoromal component functions $X$ and $r \times r$ coefficient matrix $S$}. 
Second, update $S$ by integrating 
\begin{equation} \label{eq:S_step_orth}
\partial_t S = -\left\langle X h V \right\rangle_{x,v}. 
\end{equation}
Third, let {\color{black} $L(t,v) = V(t,v) S(t)^{\top}$, a length-$r$ row vector of non-orthogonal component functions, } so that $f_{\color{black} r} = X L^{\top}$, and update $L$ by integrating 
\begin{equation} \label{eq:L_step_orth}
\partial_t L = \left\langle X h \right\rangle_x.
\end{equation}
Then orthonormalize $L$ to obtain {\color{black} an updated set of $r$ orthonoromal component functions $V$ and $r \times r$ coefficient matrix $S$.} 

%The orthogonal projection \eqref{eq:orth_proj} is the best approximation of $h$ in the tangent space of the rank-$r$ manifold $\mathcal{M}_{r}$ with respect to the $L^2$ norm. 
The scheme is memory efficient (the memory cost scales as $n^3$ instead of $n^6$, where $n$ is the number of grid points in each direction; here assumed to be equal). This can be easily seen as only low-rank factors that at most scale with $n^3$ are necessary to perform the scheme.
However, the scheme is not computationally efficient. For example to compute
$$
    \langle M(f_{\color{black} r}) V \rangle_v 
    = \frac{\rho(f_{\color{black} r})}{(2\pi T(f_{\color{black} r}))^{d_v/2}} \left\langle \exp\left(-\frac{\| v-u(f_{\color{black} r}) \|^{2}}{2T(f_{\color{black} r})}\right) V \right\rangle_v ,
$$
which is required in equation \eqref{eq:K_step_orth}, we need to, in general, compute an integral in $v$ (with cost scaling as $n^3$) for each of the $n^3$ grid points in $x$. This clearly scales as $n^6$, which is prohibitive. As explained in the introduction, various ways to overcome this issue and regain a scaling of computational cost of no more than $n^3$ have been proposed in the literature \cite{Einkemmer2019,Einkemmer2021d}. However, none of these approaches can efficiently treat the non-isothermal fully compressible case that we consider in this paper.

\section{Interpolatory dynamical low-rank approximation} 
\label{sec:interp_DLR}

In this section we propose two interpolatory dynamical low-rank schemes for the time integration of the Boltzmann-BGK equation \eqref{eq:boltzmann-bgk}. {\color{black} The key idea is to replace the orthogonal projection of $h$ onto the tangent space in \eqref{eq:orth_proj} with an oblique projection onto the same space, constructed to satisfy an interpolation property. In this section we show that} in contrast to the classic DLR algorithm outlined above, the computational cost scales as $n^3$. The first method is a direct application of the recently developed interpolatory projector-splitting scheme on low-rank matrix manifolds \cite{Dektor2024}. The second method is a variant of the rank-adaptive basis update and Galerkin scheme \cite{Ceruti2022}, where the Galerkin step is replaced by a collocation step, resulting in a new scheme we call basis update and collocate (BUC). 

\subsection{Interpolatory projector-splitting DLR} 

For the interpolatory projector-splitting integrator we impose that $f$ satisfies the Boltzmann–BGK equation \eqref{eq:boltzmann-bgk} for $r$ fixed values of $x$ and arbitrary $v$ and for $r$ fixed values of $v$ and arbitrary $x$. To achieve this, we consider 
\begin{equation} \label{eq:interp_proj} 
\partial_t f_{\color{black} r}(x,v) = 
h(x,\bm v) V(\bm v)^{-\top} V(v)^{\top} 
- X(x) X(\bm x)^{-1} h(\bm x,\bm v) V(\bm v)^{-\top} V(v)^{\top} 
+ X(x) X(\bm x)^{-1} h(\bm x, v), 
\end{equation}
where $\bm x \subset \Omega_x$ and $\bm v \subset \Omega_v$ denote (time-dependent) collections of $r$ discrete values of $x \in \Omega_x$ and $v \in \Omega_v$, respectively, and $X(\bm x)$ and $V(\bm v)$ denote the $r \times r$ matrices obtained from evaluating the basis functions at those points. 
{\color{black} The right-hand-side of \eqref{eq:interp_proj} is an oblique projection of $h$ onto a tangent space of $\mathcal{M}_r$ with the property} $\partial_t f_{\color{black} r}(\bm x,v) = h(\bm x,v)$ and $\partial_t f(x,\bm v) = h(x,\bm v)$, i.e., $\partial_t f_{\color{black} r}$ interpolates $h$ whenever $x \in \bm x$ or $v \in \bm v$. {\color{black} Hence, we refer to \eqref{eq:interp_proj} as an interpolatory projection onto the tangent space of $\mathcal{M}_r$.} Note that all functions in \eqref{eq:interp_proj} {\color{black} and the collections of discrete points $\bm x, \bm v$} are time-dependent and we suppressed the variable $t$ for convenience. One can interpret interpolatory DLR as a collocation approach on the low-rank manifold via interpolatory projection \eqref{eq:interp_proj} onto the tangent space. Meanwhile classical DLR is a Galerkin approach that minimizes the $L^2$ norm of the approximation via orthogonal projection \eqref{eq:orth_proj} onto the tangent space.  We choose points $\bm x$ and $\bm v$ at each time $t$ so that the $r \times r$ matrices $X(\bm x)$ and $V(\bm v)$ obtained from evaluating the time-dependent bases $X$ and $V$ at these points are well-conditioned. This implies that the interpolatory projection \eqref{eq:interp_proj} exists and its computation is stable. When the bases $X$ and $V$ are orthonormal, it is always possible to find a collection of points that yield well-defined interpolatory projectors using a sparse sampling algorithm. {\color{black} One can optionally increase the solution rank by selecting more than $r$ points with sparse oversampling algorithms.} We describe the selection of such points in the discrete setting in Section \ref{sec:semi-discrete}.

To integrate \eqref{eq:interp_proj} we use operator splitting, just as we have done above in \eqref{eq:K_step_orth}-\eqref{eq:L_step_orth} for the orthogonal projection, which requires solving three substeps. The interpolation points are sampled as needed within the substeps. First, ensure that the basis $V$ is orthonormal and sample $r$ points $\bm v$ so that $V(\bm v)$ is well-conditioned. Then set $K = X S$ and update $K$ by integrating 
\begin{equation} 
\partial_t K = h(x,\bm v) V(\bm v)^{-\top}.
\end{equation}
Then we perform an orthonormalization of the updated $K$ to obtain new $X$ and $S$ using a QR-decomposition. 
Second, sample $r$ points $\bm x$ from the orthonormal basis $X$ so that $X(\bm x)$ is well-conditioned. Then update $S$ by integrating  
\begin{equation}
\partial_t S(t) = - {X}(\bm x)^{-1} 
h(\bm x, \bm v) V(\bm v)^{-\top}. 
\end{equation}
Third, set $L = V S^{\top}$ and update $L$ by integrating 
\begin{equation} 
\partial_t L = h(\bm x,v)^{\top}  X(\bm x)^{-\top}. 
\end{equation}
Finally, perform an orthonormalization of the updated $L$ to obtain new $V$ and $S$ using a QR-decomposition. 

The computational advantage of the interpolatory projector-splitting integrator presented above is that solving each substep requires evaluating $h(t,x,v)$ at $r$ velocity points (K-step), $r$ spatial points (L-step) or $r$ spatial and velocity points (S-step), which is efficient even when $h$ is not expressed in a low-rank form (as is the case for the Boltzmann–BGK equation \eqref{eq:boltzmann-bgk}). 
Meanwhile the classical DLR approach that applies projector-splitting to the orthogonal tangent space projection \eqref{eq:orth_proj} requires inner products of $h$ with the bases $X$ and $V$. Such scheme requires a low-rank representation of $h$ or some other way to compute the arising integrals efficiently. 

\subsubsection{Semi-discrete interpolatory DLR scheme} \label{sec:semi-discrete} 

We discretize the spatial domain $\Omega_x$ and velocity domain $\Omega_v$ with a Fourier pseudospectral method \cite{Hesthaven_2007} using $N_x$ points $\{x^k\}_{k=1}^{N_x} \subset \Omega_x$ and and $N_v$ points $\{v^k\}_{k=1}^{N_v} \subset \Omega_v$. The discretized rank-$r$ distribution function $\bm f(t)$ is a $N_x \times N_v$ matrix with rank-$r$ representation 
\begin{equation} \label{eq:DLR_discrete}
\bm f(t) = \bm X(t) \bm S(t) \bm V(t)^{\top}, 
\end{equation}
where $\bm X(t)\in \mathbb{R}^{N_x \times r}$, $\bm S(t)\in \mathbb{R}^{r\times r}$, and $\bm V(t)\in \mathbb{R}^{N_v\times r}$. Associated with such discretization points are quadrature weights $\{w^k_x\}_{k=1}^{N_x}$ and $\{w^k_v\}_{k=1}^{N_v}$ and the discrete form of the differential operator $\nabla_x$, which is a $N_x \times N_x$ matrix that we denote by $\bm D_x$. The discrete form of \eqref{eq:h_def} is the $N_x \times N_v$ matrix 
\begin{equation} \label{eq:h_discrete}
\bm h(t) = - \left(\bm D_x \bm X(t)\right) \bm S(t) \left(\diag(\bm v) \bm V(t)\right)^{\top} + \frac{1}{\epsilon} \left(\bm M (\bm f) - \bm f \right),  
\end{equation}
with discrete Maxwellian 
\begin{equation} \label{eq:disc_maxwellian}
  \bm M(\bm f) =  \left(\frac{\bm \rho(\bm f)}{(2\pi \bm T(\bm f))^{d_v/2}} \bm 1_{N_v}^{\top} \right)\odot \exp\left(-\sum_{i=1}^{d_v} \left(\bm 1_{N_x} \bm v_i^{\top} -\bm u_i(\bm f) \bm 1_{N_v}^{\top}\right)^{2} \odot \left(\frac{1}{2\bm T(\bm f)} \bm 1_{N_v}^{\top} \right)\right), 
\end{equation}
where $\bm v_i$ is a $N_v$-dimensional column vector containing discrete values of $v_i$ (i.e.~the $i$th velocity direction), $\odot$ denotes entry-wise multiplication (Hadamard product), division of vectors is performed entry-wise, and $\bm 1_N$ denotes the $N$-dimensional column vector containing all ones. The discrete moments in \eqref{eq:disc_maxwellian} can be written using the rank-$r$ decomposition \eqref{eq:DLR_discrete} of the discrete distribution function 
\begin{equation} \label{eq:disc_moments}
\begin{aligned}
\bm \rho(\bm f) 
&= \bm f \bm w_v 
= \bm X \bm S \left(\bm V^{\top} \bm w_v\right), \\
\bm u_i(\bm f) &= \frac{1}{\bm \rho(\bm f)} \odot \left[\bm f \diag(\bm v_i) \bm w_v \right]
= \frac{1}{\bm \rho(\bm f)} \odot \left[\bm X \bm S \left(\bm V^{\top} \diag(\bm v_i) \bm w_v\right) \right], \quad i=1,\ldots, d_x, \\ 
\bm T(\bm f) &= \frac{1}{{\color{black} d_v} \bm \rho(\bm f)} \odot \left[\left( \sum_{i=1}^{d_v} \bm f \diag(\bm v_i^2) - 2\diag(\bm u_i) \bm f \diag(\bm v_i) + \diag(\bm u_i^2)\bm f \right) \bm w_v \right] \\
&= \frac{1}{{\color{black} d_v} \bm \rho(\bm f)} \odot \left[\sum_{i=1}^{d_v} \bm X \bm S\left(\bm V^{\top} \diag\left(\bm v_i^2\right) \bm w_v \right) - 2 \diag(\bm u_i) \bm X \bm S \left(\bm V^{\top} \diag(\bm v_i) \bm w_v \right) + \diag(\bm u_i^2)\bm X \bm S \left(\bm V^{\top}\bm w_v\right) \right], 
\end{aligned}
\end{equation}
where $\bm w_v$ is a column vector containing the quadrature weights $w_v^k$. 
{\color{black} In order to compare the computational cost of the interpolatory DLR and orthogonal DLR schemes, let us assume that $N_x=N_v=n$ as we have done in Section~\ref{sec:DLR}.} 
The second equality in each of the discrete moment equations \eqref{eq:disc_moments} make use of the rank-$r$ decomposition \eqref{eq:DLR_discrete} and involves matrix multiplication between matrices of dimension at most $n^3 \times r$ and Hadamard products between vectors of length $n^3$, allowing us to compute the moments with computational cost scaling as $\mathcal{O}(rn^3)$. Due to the nonlinearity in the Maxwellian, we can not compute a low-rank decomposition of \eqref{eq:h_discrete} from the rank-$r$ decomposition \eqref{eq:DLR_discrete} of $\bm f(t)$. For this reason the classical DLR approach based on the orthogonal projection \eqref{eq:orth_proj} is inefficient. The interpolatory DLR scheme does not require a low-rank form of $\bm h(t)$. Instead such scheme requires evaluating $\bm h(t)$ at a subset of $r$ row indices and $r$ column indices, which can be easily implemented with cost scaling as $\mathcal{O}(rn^3)$. 

We now describe a strategy for selecting the points $\bm x$ and $\bm v$ defining the interpolatory projection \eqref{eq:interp_proj}. In the discrete setting, sampling a point from the continuous variable $x$ corresponds to sampling an index $i \in \{1,\ldots,N_x\}$ and sampling a point from the continuous variable $v$ corresponds to sampling an index $j \in \{1,\ldots,N_v\}$. To obtain a well-conditioned tangent space projector \eqref{eq:interp_proj}, our goal is to sample $r$ indices $\I \subset \{1,\ldots,N_x\}$ and $r$ indices $\J \subset \{1,\ldots,N_v\}$ so that the condition number of the $r \times r$ matrices $\bm X\left(\I,:\right)$ and $\bm V\left(\J,:\right)$ obtained from evaluating the discrete basis functions at the indices (or continuous basis function at the discrete points) is small. For this task we employ the discrete empirical interpolation method (DEIM) \cite{Chaturantabut2010}, recalled in Algorithm \ref{alg:DEIM}, which is a greedy algorithm for minimizing the condition number of such matrices. We denote by $\texttt{DEIM}$ a subroutine that takes discretized basis functions $\bm X$ or $\bm V$ as input and returns indices $\I$ or $\J$. {\color{black} The framework proposed in this work is compatible with a variety of sparse sampling strategies, including those that offer theoretical accuracy guarantees such as \cite{cortinovis2024adaptive}. Nevertheless, for the numerical experiments presented in Section~\ref{sec:numerics} we find that the standard DEIM algorithm yields satisfactory results. It is also possible to sample more than $r$ indices using oversampling techniques such as E-DEIM \cite{hendryx2021extended} or GappyPOD+E \cite{gappy_podE}. Oversampling temporarily increases the solution rank, which can subsequently be truncated, leading to a rank-adaptive variant of the proposed projector-splitting algorithm. }

Hereafter we describe the three substeps for integrating the spatially discrete distribution function $\bm f(t)$ from time $t_0$ to time $t_1$ starting from a rank-$r$ representation $\bm f(t_0) = \bm X(t_0) \bm S(t_0) \bm V(t_0)^{\top}$. This can be considered an interpolatory projector splitting integrator. 

\begin{itemize}

\item {\bf K-step}: update $\bm X(t_0) \to \bm X(t_1)$ and $\bm S(t_0) \to \bm R(t_1)$. 

\vs
Compute interpolation indices $\J = \DEIM(\bm V(t_0))$. Then integrate the $N_x \times r$ differential equation 
\begin{equation} \label{eq:Kstep} 
\frac{d \bm K(t)}{dt} = 
\bm h\left(t,:,\J\right) \left[\bm V\left(t_0,\J,:\right)\right]^{-\top}, \qquad \bm K(t_0) = \bm X(t_0)\bm S(t_0), 
\end{equation}
from $t_0$ to $t_1$, and perform a QR-decomposition $\bm K(t_1) = \bm X(t_1) \bm R(t_1)$. 

\item {\bf S-step}: update $\bm R(t_1) \to \tilde{\bm S}(t_1)$. 

\vs
Compute interpolation indices $\I = \DEIM(\bm X(t_1))$. Then integrate the $r \times r$ matrix differential equation 
\begin{equation} \label{eq:Sstep}
\frac{d\tilde{\bm S}(t)}{dt} = -\left[\bm X\left(t_1,\I,:\right)\right]^{-1} 
\bm h\left(t,\I,\J\right) \left[\bm V\left(t_0,\J,:\right)\right]^{-\top}, 
\qquad \tilde{\bm S}(t_0) = \bm R(t_1), 
\end{equation}
from $t_0$ to $t_1$. 

\item {\bf L-step}: update $\bm V(t_0) \to \bm V(t_1)$ and $\tilde{\bm S}(t_1) \to \bm S(t_1)$. 

\vs
Integrate the $N_v \times r$ matrix differential equation 
\begin{equation} \label{eq:Lstep}
\frac{d \bm L(t)}{dt} = 
\bm h\left(t,\I,:\right)^{\top} 
\left[\bm X\left(t_1,\I,:\right)\right]^{-\top}, \qquad \bm L(t_0) = \bm V(t_0)\tilde{\bm S}(t_1)^{\top}, 
\end{equation}
from time $t_0$ to $t_1$, and perform a QR-decomposition $\bm L(t_1) = \bm V(t_1) \bm S(t_1)^{\top}$. 

\end{itemize}

\noindent
The rank-$r$ distribution function at time $t_1$ is $\bm f(t_1) = \bm X(t_1) \bm S(t_1) \bm V(t_1)^{\top}$. 

In the $\bm K$-step the computational cost for computing $\J$ using $\texttt{DEIM}$ is $\mathcal{O}(rn^3)$. Integrating the differential equation \eqref{eq:Kstep}, e.g., using an Euler forward scheme requires evaluating $\bm h(t,:,\J)$ once which costs $\mathcal{O}(rn^3)$, inverting a $r \times r$ matrix $\mathcal{O}(r^3)$, and updating $\bm K$ which costs $\mathcal{O}(rn^3)$. The cost of the QR decomposition to obtain $\bm X(t_1)$ and $\bm R(t_1)$ is $\mathcal{O}(r^2n^3)$. Thus the cost for the $\bm K$-step scales as $\mathcal{O}(r^2n^3)$. In the $\bm S$-step the cost of computing $\I$ with $\texttt{DEIM}$ is $\mathcal{O}(rn^3)$. Integrating the differential equation \eqref{eq:Sstep} with Euler forward requires evaluating $\bm h(t,\I,\J)$, inverting two $r \times r$ matrices $\mathcal{O}(r^3)$, and updating $\bm S$ costing $\mathcal{O}(r^2)$. The total cost of the $\bm S$-step scales as $\mathcal{O}(rn^3 + r^3)$. In the $\bm L$-step, integrating \eqref{eq:Lstep} with Euler forward requires evaluating $\bm h(t,\I,:)$ costing $\mathcal{O}(rn^3)$, inverting a $r \times r$ matrix $\mathcal{O}(r^3)$, and updating $\bm L$ which costs $\mathcal{O}(rn^3)$. Then the QR decompsition to obtain $\bm V(t_1)$ and $\bm S(t_1)$ costs $\mathcal{O}(r^2n^3)$. The total cost of the $\bm L$-step is $\mathcal{O}(r^2n^3)$. The total cost of one step of the interpolatory projector-splitting method for the Boltzmann-BGK equation \eqref{eq:boltzmann-bgk} is $\mathcal{O}(r^2n^3 + r^3)$. 

\begin{algorithm}
 \caption{DEIM index selection (adapted from \cite{Sorensen_2016}).}
\label{alg:DEIM}
\begin{algorithmic}[1]
\Require 
%\Statex
    $\bm M \in \mathbb{R}^{n \times r}$ with $n \geq r$
\Ensure 
	$\bm l=(l_1,\ldots,l_r) \subset \{1,\ldots,n\}$ minimizing condition number of $\bm M(\bm l,:)$ 
%    \Statex 
%
%
\State $\bm m = \bm M(:,1)$
\State $l_1 = \arg\max(\text{abs}(\bm m))$
\For{$j = 2,3,\ldots,r$ }
        \State $\bm m = \bm M(:,j)$ 
        \State $\bm c = \bm M(l,1:j-1)^{-1} \bm m(\bm l)$
        \State $\bm r =\bm m - \bm M(:,1:j-1) \bm c$ 
        \State $l_j = \arg \max(\text{abs}(\bm r))$
\EndFor
\end{algorithmic}
\end{algorithm}

\subsection{The basis update and collocate (BUC) integrator} \label{sec:BUC}

As an alternative to the orthogonal projector-splitting DLR summarized in Section \ref{sec:DLR}, an unconventional scheme was introduced in \cite{Ceruti2022}. The idea of this scheme is to update the low-rank factors $X$ and $V$ and then perform a Galerkin projection onto the updated bases. Thus the scheme is referred to as the basis update and Galerkin (BUG) integrator. The integrator in \cite{Ceruti2022}, however, needs to project the initial condition onto the new basis. This introduces numerical error that in some situations can be detrimental to the accuracy and stability of the numerical scheme, see, e.g., \cite{Einkemmer2023}. This problem was remedied in \cite{Ceruti2022a} by using a larger (i.e.~augmented) basis that uses both the initial as well as the newly computed basis functions. 
{\color{black} The BUC integrator proposed in this section follows the same steps as the augmented BUG integrator \cite{Ceruti2022a}, which we recall hereafter for the convenience of the reader. 
\begin{itemize}

\item {\bf K-step}: update $\bm X(t_0) \to \bm X(t_1)$. 

\vs
Integrate the $N_x \times r$ differential equation 
\begin{equation} \label{eq:BUG_K}
\frac{d \bm K(t)}{dt} =
\bm h \bm V\left(t_0\right), \qquad \bm K(t_0) = \bm X(t_0)\bm S(t_0), 
\end{equation}
from time $t_0$ to $t_1$. Then perform a QR-decomposition of the $N_x \times 2r$ augmented matrix $\left[ \bm K(t_1), \bm X(t_0) \right]= \widehat{\bm X}(t_1) \bm R$ and compute the $2r \times r$ matrix $\widehat{\bm M} = \widehat{\bm X}(t_1)^{\top} \bm X(t_0)$. The $\bm R$ is discarded.

\item {\bf L-step}: update $\bm V(t_0) \to \bm V(t_1)$. 

\vs
Integrate the $N_v \times r$ matrix differential equation 
\begin{equation} \label{eq:BUG_L}
\frac{d \bm L(t)}{dt} = \bm h^{\top} \bm X(t_0), \qquad \bm L(t_0) = \bm V(t_0)\bm S(t_0)^{\top}, 
\end{equation}
from time $t_0$ to $t_1$. Then perform a QR-decomposition of the $N_v \times 2r$ augmented matrix $\left[\bm L(t_1), \bm V(t_0) \right] = \widehat{\bm V}(t_1) \tilde{\bm R}$ and compute the $2r \times r$ matrix $\widehat{\bm N} = \widehat{\bm V}(t_1)^{\top} \bm V(t_0)$. The $\tilde{\bm R}$ is again discarded. 

\item {\bf S-step}: update $\bm S(t_0) \to \widehat{\bm S}(t_1)$. 

\vs 
Integrate the $2r \times 2r$ matrix differential equation 
\begin{equation} \label{eq:BUG_S}
\frac{d \widehat{\bm S}(t)}{dt} 
= \widehat{\bm X}^{\top}
\bm h 
\widehat{\bm V}, 
\qquad \widehat{\bm S}(t_0) = \widehat{\bm M} \bm S(t_0) \widehat{\bm N}^{\top},
\end{equation}
from time $t_0$ to $t_1$. 

\item {\bf Truncate}:

\vs
Compute the SVD $\widehat{\bm S}(t_1) = \widehat{\bm P}\widehat{\bm \Sigma} \widehat{\bm Q}^{\top}$ and truncate singular values with relative tolerance $\vartheta$ to obtain a new rank $r_1 \leq 2r$. The new factors for the rank-$r_1$ approximation of $\bm f(t_1)$ are obtained as follows. Let $\bm S(t_1)$ be the diagonal matrix with the first $r_1$ singular values of $\widehat{\bm S}(t_1)$. Let $\bm P= \widehat{\bm P}(:,1:r_1)$ and $\bm Q = \widehat{\bm Q}(:,1:r_1)$ be the first $r_1$ columns of $\widehat{\bm P}$ and $\widehat{\bm Q}$, respectively. 
Finally let $\bm X(t_1) = \widehat{\bm X}(t_1) \bm P \in \mathbb{R}^{N_x \times r_1}$ and $\bm V(t_1) = \widehat{\bm V}(t_1)\bm Q \in \mathbb{R}^{N_v \times r_1}$. 

\end{itemize}
}

The augmented BUG scheme has a number of advantages over the projector-splitting integrator. The augmented bases allow for a natural adaptive rank increase of the solution during each time step. Such rank-adaptivity is useful to ensure an accurate approximation during time integration, although we note that rank-adaptive variants of the projector splitting integrator have been developed \cite{Hochbruck2023,Dektor2021}. Second, the $\bm S$-step in the projector-splitting schemes (both orthogonal and interpolatory) has a minus sign making it a ``backwards'' step in time. This can be unstable for dissipative problems. The BUG integrator does not have a backwards step. Finally, the $\bm K$ and $\bm L$ steps of the BUG scheme can be performed in parallel. We note, however, that the larger approximation space increases computational cost compared to the projector splitting (and also the classic BUG) integrator. We will study this in more detail in section \ref{sec:numerics}.

Hereafter we propose a rank-adaptive interpolatory DLR scheme inspired by the augmented BUG DLR scheme \eqref{eq:BUG_K}-\eqref{eq:BUG_S}. The {\color{black} steps are almost identical, but the orthogonal projections onto the $\bm X$ and $\bm V$ bases are replaced with interpolatory projections, which allows us to efficiently perform the updates in each substep for nonlinear equations such as the Boltzmann-BGK equation \eqref{eq:boltzmann-bgk}. Below are the steps for our BUG-inspired interpolatory scheme. } 

\begin{itemize}

\item {\bf K-step}: update $\bm X(t_0) \to \bm X(t_1)$. 

\vs
Compute interpolation indices $\J = \DEIM(\bm V(t_0))$ and integrate the $N_x \times r$ differential equation 
\begin{equation}
\frac{d \bm K(t)}{dt} =
\bm h\left(t,:,\J\right) \bm V\left(t_0,\J,:\right)^{-\top}, \qquad \bm K(t_0) = \bm X(t_0)\bm S(t_0), 
\end{equation}
from time $t_0$ to $t_1$. Then perform a QR-decomposition of the $N_x \times 2r$ augmented matrix $\left[ \bm K(t_1), \bm X(t_0) \right]= \widehat{\bm X}(t_1) \bm R$ and compute the $2r \times r$ matrix $\widehat{\bm M} = \widehat{\bm X}(t_1)^{\top} \bm X(t_0)$. The $\bm R$ is discarded.

\item {\bf L-step}: update $\bm V(t_0) \to \bm V(t_1)$. 

\vs
Compute interpolation indices $\I = \DEIM(\bm X(t_0))$ and 
integrate the $N_v \times r$ matrix differential equation 
\begin{equation}
\frac{d \bm L(t)}{dt} = \bm h(t,\I,:)^{\top} \bm X(t_0,\I,:)^{-\top}, \qquad \bm L(t_0) = \bm V(t_0)\bm S(t_0)^{\top}, 
\end{equation}
from time $t_0$ to $t_1$. Then perform a QR-decomposition of the $N_v \times 2r$ augmented matrix $\left[\bm L(t_1), \bm V(t_0) \right] = \widehat{\bm V}(t_1) \tilde{\bm R}$ and compute the $2r \times r$ matrix $\widehat{\bm N} = \widehat{\bm V}(t_1)^{\top} \bm V(t_0)$. The $\tilde{\bm R}$ is again discarded. 

\item {\bf S-step}: update $\bm S(t_0) \to \widehat{\bm S}(t_1)$. 

\vs 
Compute $2r$ indices $\widehat{\I} = \DEIM\left(\widehat{\bm X}(t_1) \right)$ and $2r$ indices $\widehat{\J} = \DEIM\left(\widehat{\bm V}(t_1)\right)$ for interpolatory projections onto the augmented bases. Then integrate the $2r \times 2r$ matrix differential equation 
\begin{equation}
\frac{d \widehat{\bm S}(t)}{dt} 
= \widehat{\bm X}\left(t_1,\widehat{\I},:\right)^{-1}
\bm h\left(t,\widehat{\I},\widehat{\J}\right) 
\widehat{\bm V}\left(t_1,\widehat{\J},:\right)^{-\top}, 
\qquad \widehat{\bm S}(t_0) = \widehat{\bm M} \bm S(t_0) \widehat{\bm N}^{\top},
\end{equation}
from time $t_0$ to $t_1$. 

\item {\bf Truncate}:

\vs
Compute the SVD $\widehat{\bm S}(t_1) = \widehat{\bm P}\widehat{\bm \Sigma} \widehat{\bm Q}^{\top}$ and truncate singular values with relative tolerance $\vartheta$ to obtain a new rank $r_1 \leq 2r$. The new factors for the rank-$r_1$ approximation of $\bm f(t_1)$ are obtained as follows. Let $\bm S(t_1)$ be the diagonal matrix with the first $r_1$ singular values of $\widehat{\bm S}(t_1)$. Let $\bm P= \widehat{\bm P}(:,1:r_1)$ and $\bm Q = \widehat{\bm Q}(:,1:r_1)$ be the first $r_1$ columns of $\widehat{\bm P}$ and $\widehat{\bm Q}$, respectively. 
Finally let $\bm X(t_1) = \widehat{\bm X}(t_1) \bm P \in \mathbb{R}^{N_x \times r_1}$ and $\bm V(t_1) = \widehat{\bm V}(t_1)\bm Q \in \mathbb{R}^{N_v \times r_1}$. 

\end{itemize}

\noindent
The rank-$r_1$ distribution function at time $t_1$ is $\bm f(t_1) = \bm X(t_1) \bm S(t_1) \bm V(t_1)^{\top}$. 

{\color{black} The $\bm K$-step and $\bm L$-step construct orthogonal (augmented) bases $\widehat{\bm X}(t_1)$ and $\widehat{\bm V}(t_1)$, respectively, just as in the BUG integrator described above. However, the evolution equations used to update the bases require only a subset of entries from $\bm h$, making these basis update steps efficient for nonlinear equations for which the BUG integrator is impractical. } Following the BUG integrator, we augment with $\bm X(t_0),\bm V(t_0)$ so that the initial value $\bm f(t_0)$ can be represented exactly. In the $\bm S$-step we compute a coefficient matrix $\widehat{\bm S}(t)$ so that $\partial_t \bm f(t)$ interpolates $\bm h(t)$ whenever $i \in \widehat{\I}$ or $j \in \widehat{\J}$ for all $t \in [t_0,t_1]$. Hence we refer to the proposed scheme as basis update and collocate (BUC), which shares the same advantageous properties as the BUG scheme. {\color{black} Namely the $\bm K$-step and $\bm L$-step can be performed in parallel, there is no backwards in time step, and the scheme is naturally rank-adaptive. }
Note that it may be advantageous to augment the bases in the $\bm K,\bm L$ steps using columns or rows from the right-hand side $\bm h$ in addition to the initial values $\bm X(t_0),\bm V(t_0)$ so that the indices $\widehat{\I},\widehat{\J}$ are chosen with information from $\bm h$. For our numerical demonstrations we do not perform this additional augmentation and obtain good results with the scheme above. Finally the truncation step allows us to control the accuracy of the low-rank approximation by adaptively selecting the rank of the distribution function at each time step. 
{\color{black} The computational complexity of the BUC scheme is the same as the interpolatory projector-splitting DLR scheme $\mathcal{O}(r^2n^3+r^3)$ for the Boltzmann-BGK equation \eqref{eq:boltzmann-bgk}. However, due to the temporary doubling of ranks, the computational cost involves a larger constant factor. The comparison of computational cost for the interpolatory projector-splitting and BUC schemes is discussed further in Section~\ref{sec:shear_2d2v}. }

\section{Numerical Results} \label{sec:numerics}

\begin{figure}[!t]
\centering
\includegraphics[scale=0.28]{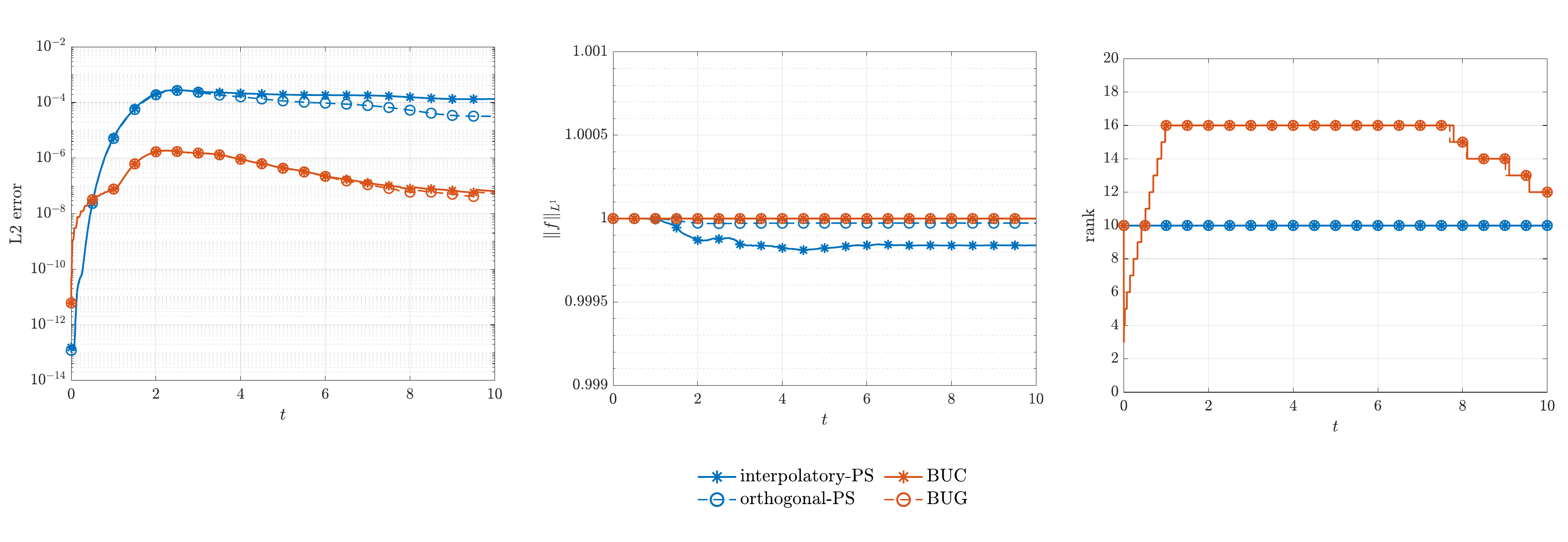}
\caption{\color{black} A comparison of the proposed interpolatory DLR schemes with their corresponding orthogonal versions for the 1d1v toy problem. Shown is the $L^2$ error, rank, and mass versus time for the interpolatory projector-splitting, Galerkin projector-splitting, BUC, and BUG schemes.} 
\label{fig:1d1v_error}
\end{figure}

\begin{figure}[!t]
\centering
\includegraphics[scale=0.13]{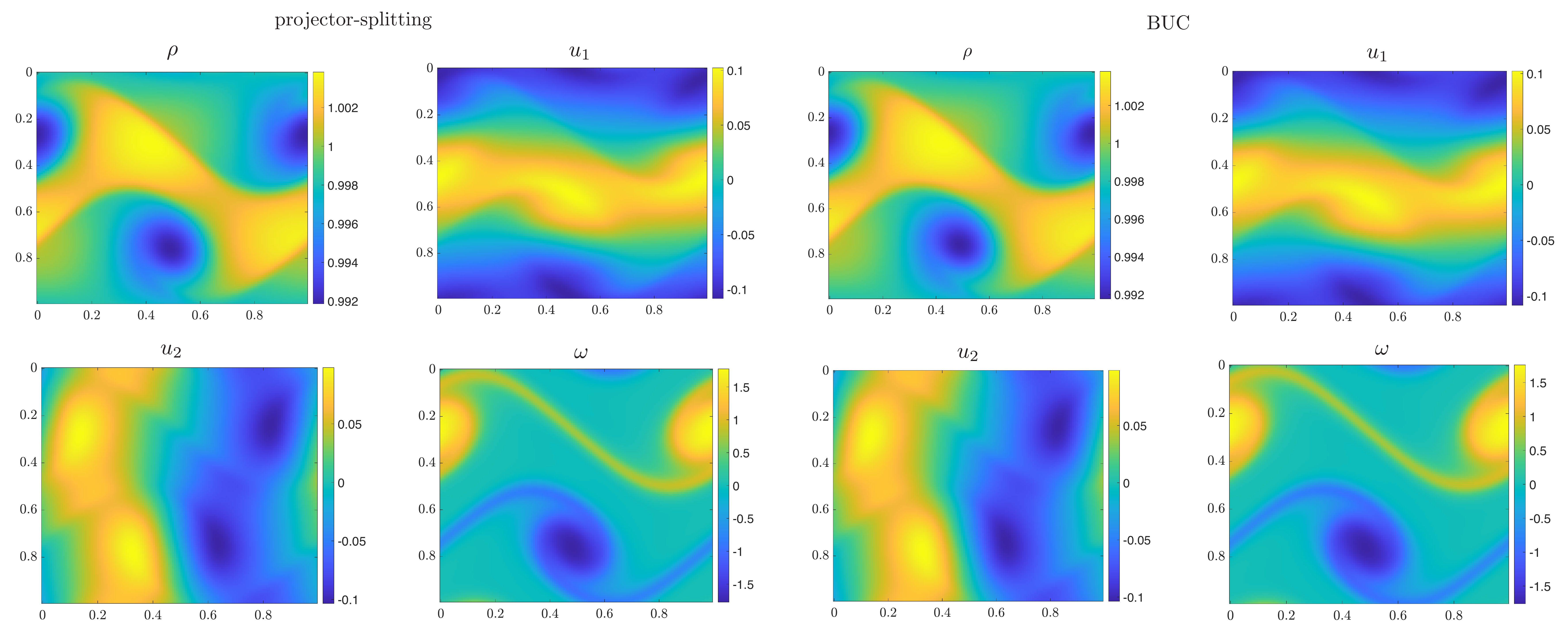}
\caption{Density, velocities, and vorticity of the 2d2v shear flow with $\epsilon=10^{-4}$ at time $t=12$ computed using the low-rank interpolatory projector-splitting integrator (left) and BUC integrator (right). We used $128$ points in each spatial direction, $16$ points in each velocity direction, and time step-size $\Delta t = 10^{-4}$. }
\label{fig:2d2v_shear_BUC}
\end{figure}

\begin{figure}[!t]
\centering
\includegraphics[scale=0.35]{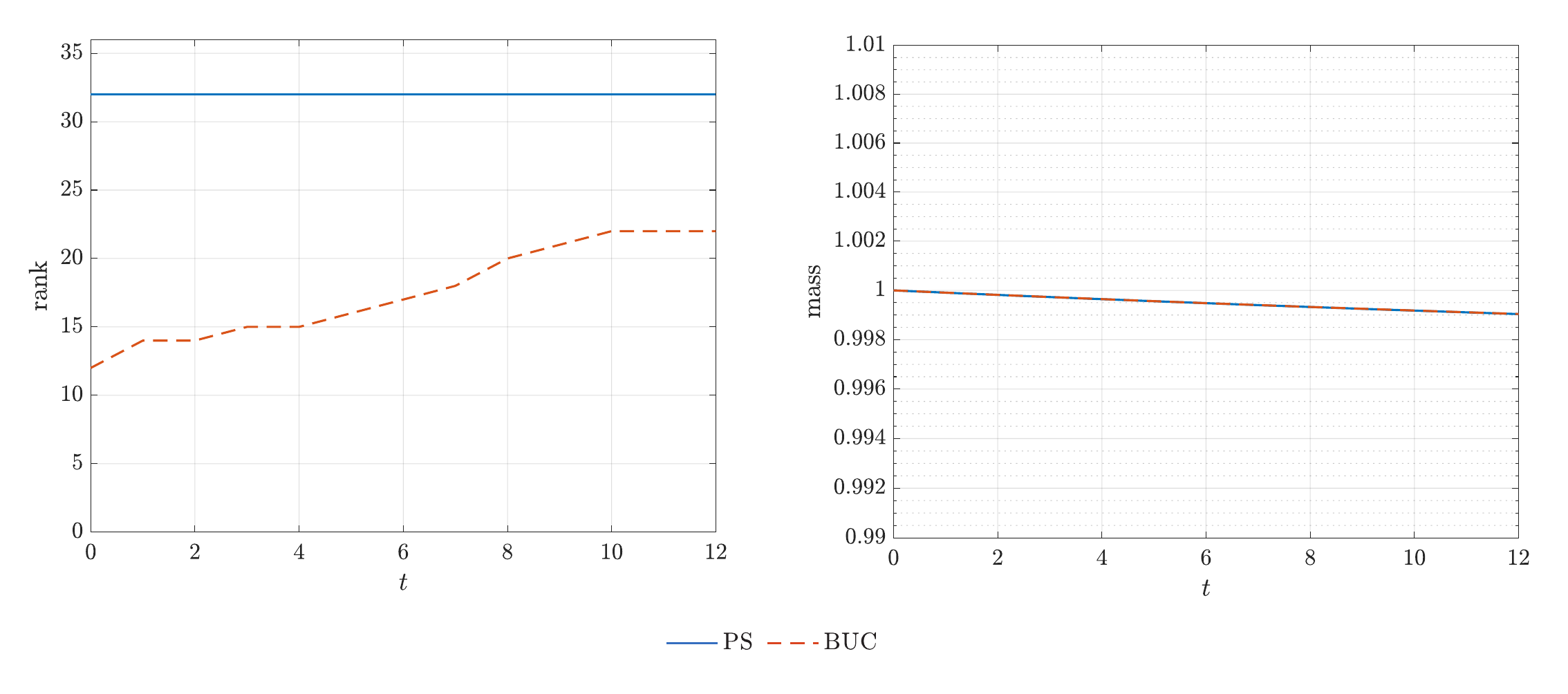}
\caption{Rank and mass versus time for the 2d2v shear flow with $\epsilon=10^{-4}$ computed with the BUC and projector-splitting schemes. The BUC scheme is rank-adaptive with singular value threshold $\vartheta=10^{-6}$, while the projector-splitting scheme uses a fixed-rank. Both methods lose mass on the order of $10^{-3}$. } 
\label{fig:2d2v_rank_mass}
\end{figure}

In this section we provide several examples in the {\color{black} 1d1v, } 2d2v, and 3d3v settings demonstrating the proposed low-rank methods for the Boltzmann-BGK equation. All numerical experiments were performed on a Macbook Pro with M2 chip and 64GB RAM. 

{\color{black}
\subsection{1d1v toy problem}

We begin by validating the proposed interpolatory DLR schemes and comparing them with the corresponding orthogonal DLR schemes on a 1d1v example. We consider the initial distribution 
\begin{equation} 
\label{eq:1d1v_ic}
f(0,x,v) = \frac{1}{m}\left[ 1 + \exp\left(-(x -x_0)^2/(2)\right) \right]
\frac{1}{\sqrt{2\pi}} \exp\left(-v^2/2\right),
\end{equation}
where $m$ is a normalization constant so that $\|f(0,x,v)\|_{L^1}=1$. We discretize the spatial domain $\Omega_x = [-6,6]$ and velocity domain $\Omega_v = [-6,6]$ using $128$ pseudospectral points in each domain and set $\epsilon=1$ in the Boltzmann-BGK equation. In this case the solution has few enough degrees of freedom to be stored as a full matrix at each time step. Thus we compute a reference solution by integrating a full matrix solution with an explicit RK4 time integration scheme with step-size $\Delta t = 10^{-3}$. 

We also compute dynamical low-rank solutions to compare the proposed interpolatory DLR schemes with the corresponding DLR schemes using orthogonal projections. Since the discrete distribution function $\bm f(t)$ and its time-derivative $d \bm f(t)/dt = \bm h(t)$ can be stored as full matrices, in this setting we can compute the orthogonal projections onto low-rank manifolds required for such DLR schemes allowing us to directly compare the interpolatory and orthogonal DLR schemes. However, for the larger problems considered below such projections are prohibitively expensive. The initial distribution function \eqref{eq:1d1v_ic} is rank-$1$. Each DLR scheme is initialized with a rank-$10$ representation of the initial condition using randomly generated orthogonal component functions with zero coefficients. The interpolatory projector-splitting and orthogonal projector-splitting integrators use constant rank $10$ for all time. Meanwhile the BUC and BUG schemes are rank-adaptive using a singular value threshold of $\vartheta = 10^{-6}$ at each time step and maximum rank of $16$. All substeps of the DLR schemes are integrated with explicit RK4 using step-size $\Delta t = 10^{-3}$. 

In Figure~\ref{fig:1d1v_error} (left) we plot the $L^2$ error of each DLR solution when compared with the full matrix reference solution. We observe that the errors of the interpolatory schemes agree with the corresponding orthogonal schemes until around $t=4$ after which the error in the interpolatory-PS solution becomes slightly larger than the orthogonal-PS solution. This is expected since the orthogonal DLR schemes compute the best approximation of $d\bm f(t)/dt$ in the tangent space each time in the $L^2$ norm, whereas the interpolatory schemes compute a sub-optimal projection. A similar deviation in error is observed amongst the BUC and BUG schemes around time $t=7$. 
In Figure~\ref{fig:1d1v_error} (middle) we plot the mass versus time for each DLR solution. The interpolatory-PS solution loses the most mass, but the error is still on the order of $10^{-3}$, which is reasonable. 
In Figure~\ref{fig:1d1v_error} (right) we plot the ranks of the DLR solutions. As we mentioned above, the projector-splitting schemes are constant rank while the BUC/BUG schemes are rank-adaptive with singular value tolerance $\vartheta = 10^{-6}$. We observe that the solutions ranks selected by the BUC and BUG schemes are almost identical at each time step. The BUC scheme takes a few additional time steps to reduce the rank after time $t=7$ as the distribution function approaches equilibrium. Overall, we observe the accuracy of the interpolatory schemes is only slightly worse than the orthogonal DLR schemes.

}

\subsection{Shear flow}

{\color{black} Next we consider} shear flows in 2d2v and 3d3v, a classic example in fluid dynamics. 

\subsubsection{2d2v} \label{sec:shear_2d2v}
We simulate the 2d2v shear flow considered in \cite[Section 7.1]{Einkemmer2021d}. The initial values are 
\begin{equation}
\begin{aligned}
\rho(0,x_1,x_2) = 1, \quad 
u_1(0,x_1,x_2) = v_0 \begin{cases}
\tanh\left(\frac{x_2-1/4}{\Delta} \right), \quad x_2 \leq \frac{1}{2}, \\
\tanh\left(\frac{3/4 - x_2}{\Delta} \right), \quad x_2 > \frac{1}{2}, 
\end{cases}
\quad u_2(0,x_1,x_2) = \delta \sin(2\pi x_1),
\end{aligned}
\end{equation}
where $(x_1,x_2) \in [0,1]^2$, $v_0=0.1$, $\Delta=1/30$ and $\delta = 5 \cdot 10^{-3}$. Following \cite{Einkemmer2021d} we set $\epsilon=10^{-4}$, which corresponds to a Reynolds number $\text{Re}=1000$. We computed low-rank solutions using the interpolatory projector-splitting scheme. This algorithm is not rank-adaptive and we thus choose a fixed rank of $32$ for all time $t$. We also computed a low-rank solution using the BUC scheme, which is rank-adaptive. We used an initial rank of $12$ and set the singular value truncation threshold to $\vartheta=10^{-6}$ for each time step. In Figure~\ref{fig:2d2v_rank_mass} we plot the rank and mass versus time for the projector-splitting and BUC solutions. The adaptively selected rank of the BUC solution remains smaller than the rank of the projector-splitting solution for all time. The total CPU-time for integrating the projector-splitting solution to time $t=12$ was approximately $7.8$ hours while the total CPU-time for the BUC scheme was approximately $6.6$ hours. Although a single time-step of the projector-splitting integrator is less expensive than a single step of the BUC integrator for solutions of the same rank, the rank-adaptive BUC integrator is faster overall since the rank of the solution is smaller for all time. To compare the cost of the interpolatory methods on equal footing, we also ran the BUC solver with fixed rank $32$ for all time $t$ which took approximately 11.1 hours. The increased computational cost of the BUC scheme compared to the interpolatory projector-splitting scheme is due to the doubling of the ranks, which is truncated after the {\bf S-step}. We observe that both solutions are not mass conservative, but only make an error on the order of $10^{-3}$ while integrating to $t=12$.

In Figure \ref{fig:2d2v_shear_BUC} we plot the density $\rho$, velocities $u_1,u_2$ and vorticity $\omega = \partial_{x_1} u_2 - \partial_{x_2} u_1$ at time $t=12$ from the solutions obtained with the projector-splitting and BUC schemes. For both schemes we observe excellent agreement of the velocities and vorticity with the plots shown in \cite{Einkemmer2021d}. For the density $\rho$, we observe a good agreement with slight deviation. 

\subsubsection{3d3v}

\begin{figure}[!t]
\centering
\includegraphics[scale=0.4]{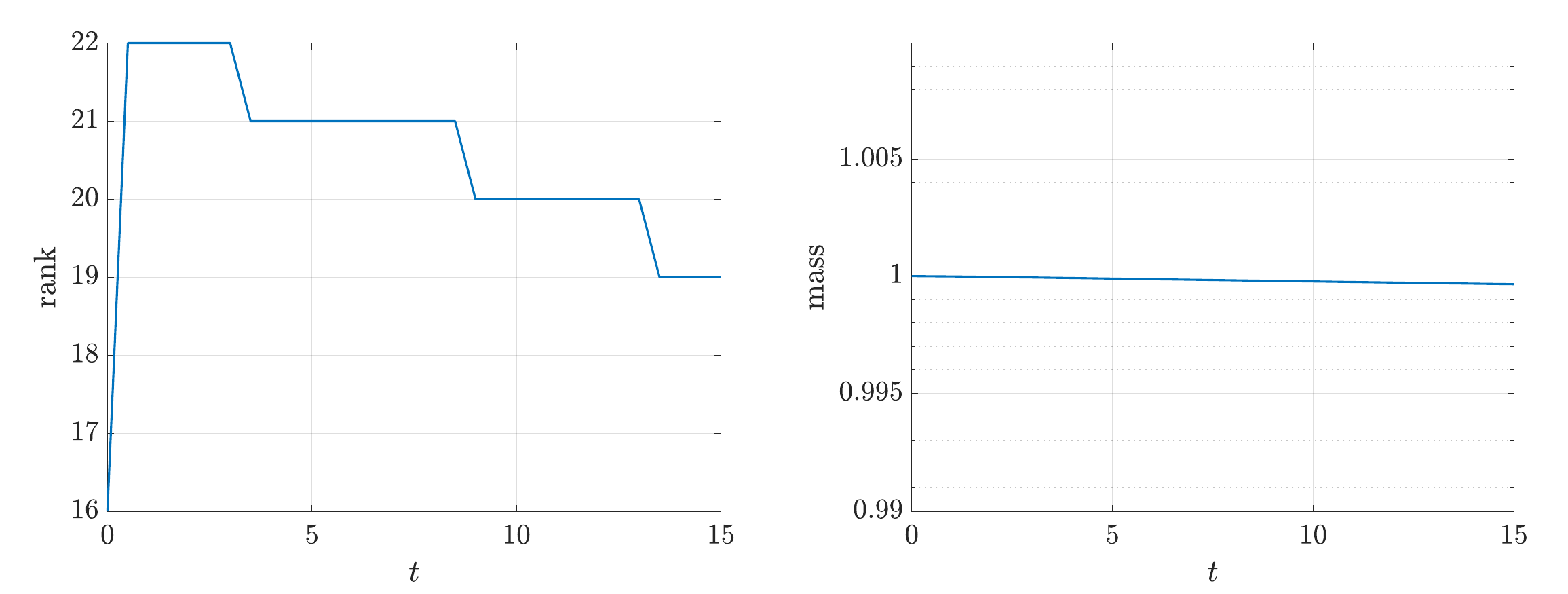}
\caption{Rank and mass versus time for the 3d3v shear flow with $\epsilon=5 \cdot 10^{-4}$ computed with the BUC scheme using rank-adaptive threshold $\vartheta=10^{-5}$. } 
\label{fig:3d3v_shear_rank_mass}
\end{figure}

Next we consider a shear flow in the 3d3v setting. For the initial values we define 
\begin{equation}
\phi(r) = v_0\tanh\left(\frac{1/4 - r^2}{\Delta} \right)
\end{equation}
and then set $\rho=1$, $T=1$, and
\begin{equation} \label{eq:vel_3d3v_shear}
\begin{aligned}
u_1(0,x_1,x_2,x_3) &= \phi(r), \\
u_2(0,x_1,x_2,x_3) &= \delta \sin(2\pi x_1) \cos(\theta), \\
u_3(0,x_1,x_2,x_3) &= \delta \sin(2\pi x_1) \sin(\theta), 
\end{aligned}
\end{equation}
where $r=\sqrt{x_2^2+x_3^2}$ and $\theta$ is the corresponding angle so that $(r,\theta)$ form polar coordinates in the $(x_2,x_3)$ plane. As computational domain we take $(x_1, x_2, x_3) \in \Omega_x = [0,1] \times [-1,1] \times [-1,1]$. The velocities \eqref{eq:vel_3d3v_shear} separate $\Omega_x$ into two regions. For roughly $r = \sqrt{x_2^2+x_3^2}\leq 1/2$ the fluid moves with $u_3=1$, while outside of that cylinder the fluid moves with $u_3=-1$. Since $(\cos \theta, \sin \theta)$ points along the normal of the interface between the two regions, the perturbation is orthogonal to the direction of the main flow and varies in the direction of the flow (as for the 2d2v shear flow above). 

For our numerical demonstration we set parameters in the initial values as $\Delta = 1/30$, $\delta = 10^{-3}$ and $v_0 = 0.1$. We discretized the spatial domain $\Omega_x$ with a tensor product grid using $100$ points in each dimension and the velocity domain $\Omega_v = [-6,6]^3$ with a tensor product grid using $16$ points per dimension. Thus the total number of degrees of freedom in a single time snapshot of the spatially discrete distribution function $\bm f$ without low-rank compression is $100^3 \cdot 16^3 \approx 4 \cdot 10^9$ requiring approximately $32$ gigabytes of memory. To avoid computing the SVD of the initial discretized distribution function $\bm f_0$ we obtain a low-rank approximation using the TT-cross maximum volume algorithm~\cite{TT-cross} with error tolerance $10^{-10}$. This yields a rank $16$ approximation of the initial discrete distribution function requiring approximately $130$ megabytes of memory (a compression rate of approximately 1:250). 

\begin{figure}[!t]
\centering
\includegraphics[scale=0.26]{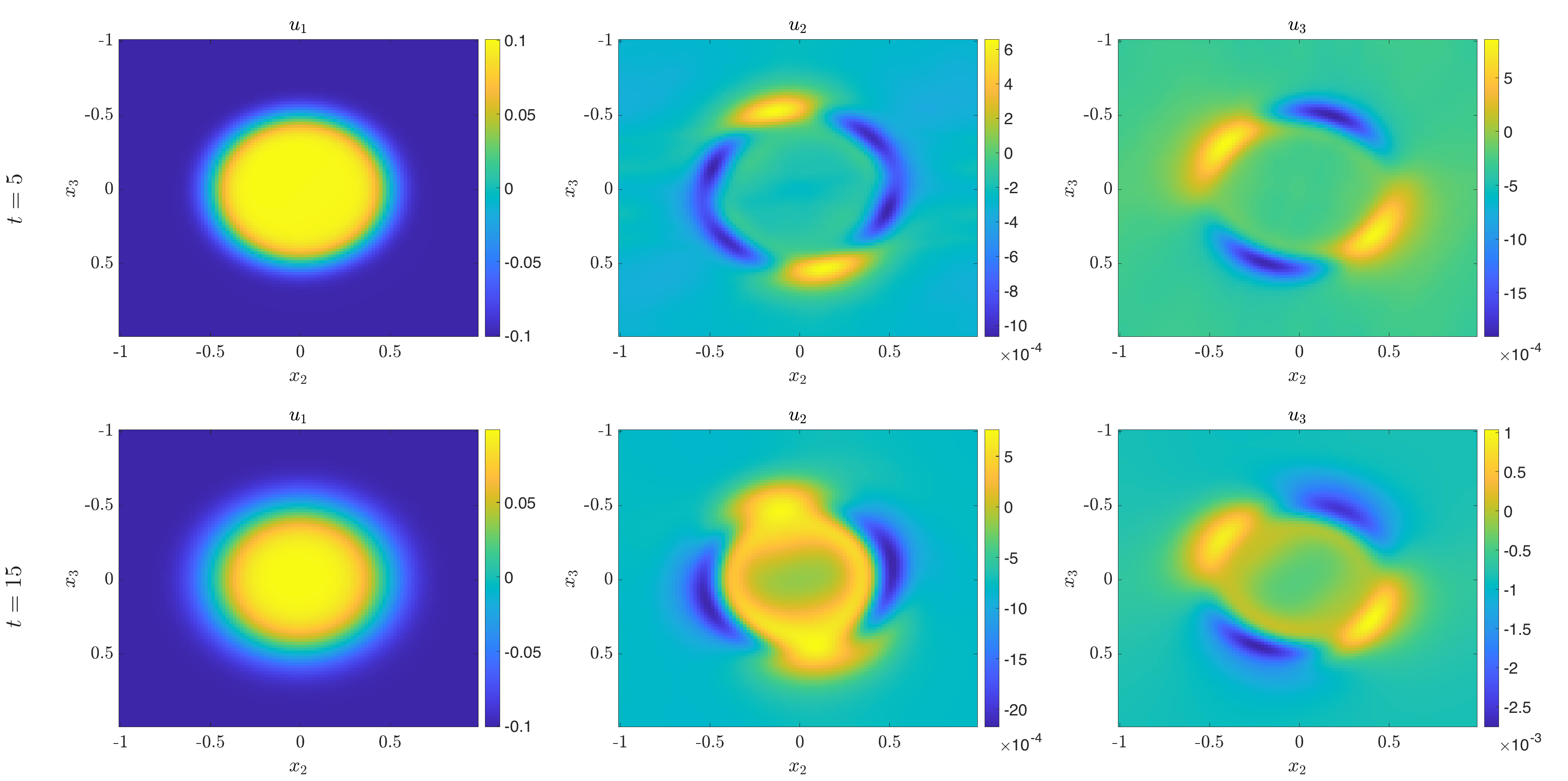}
\caption{$(x_2,x_3)$-slices of the velocity components for the 3d3v shear flow solution computed using the rank-adaptive BUC scheme. Results are shown at $x_1\approx 0.49$ for $t=5$ (top row) and $t=15$ (bottom row). }
\label{fig:3d3v_shear_vel}
\end{figure}

\begin{figure}[!t]
\centering
\includegraphics[scale=0.26]{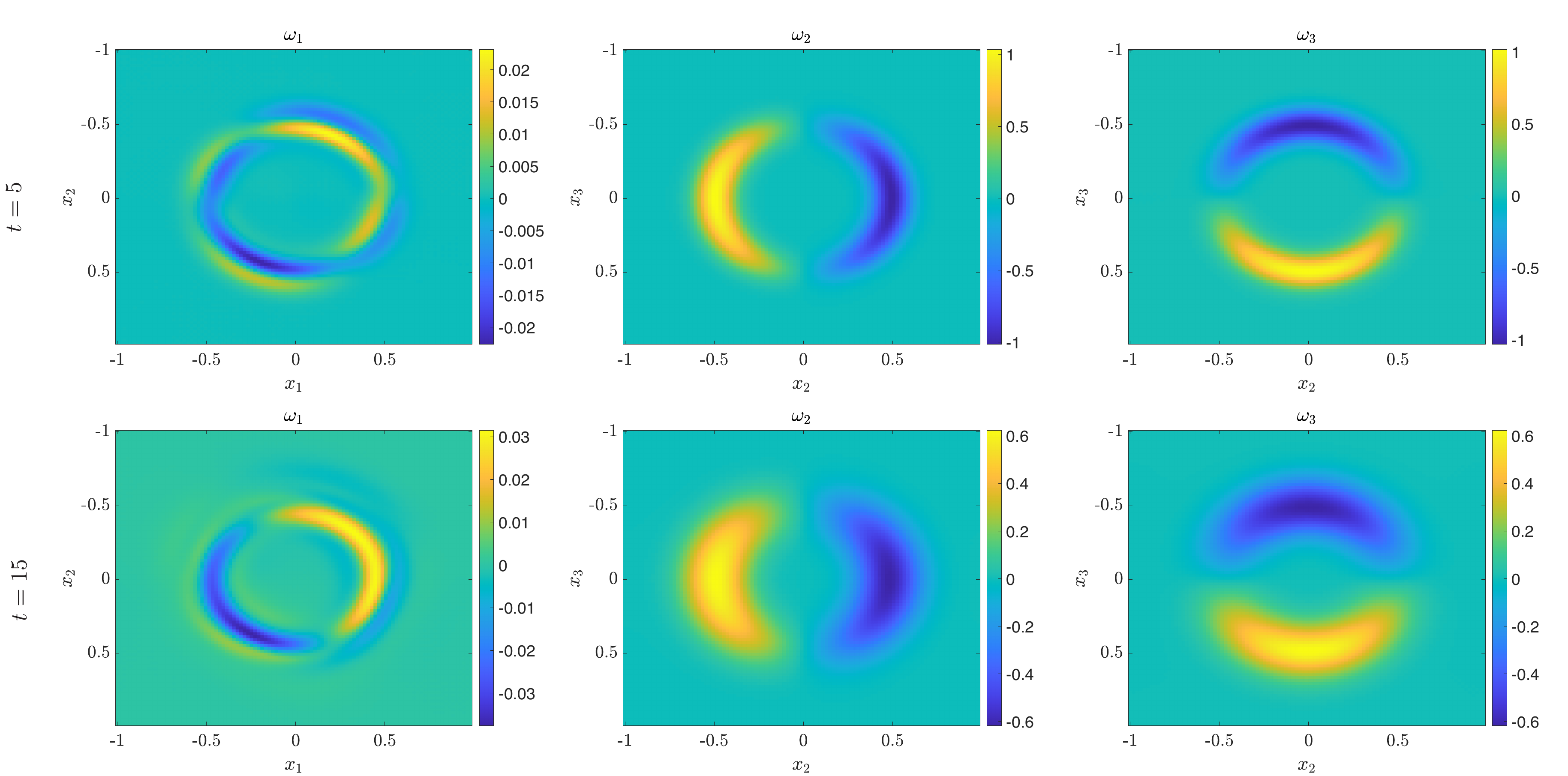}
\caption{$(x_2,x_3)$-slices of the vorticity components for the 3d3v shear flow solution computed using the rank-adaptive BUC scheme. Results are shown at $x_1\approx 0.49$ for $t=5$ (top row) and $t=15$ (bottom row). }
\label{fig:3d3v_shear_vor}
\end{figure}

We set $\epsilon= 5\cdot 10^{-4}$, corresponding to a Reynolds number of $200$, in the Boltzmann-BGK equation \eqref{eq:boltzmann-bgk} and integrate from $t=0$ to $t=15$ using the BUC scheme described in Section~\ref{sec:BUC} with step-size $\Delta t = 10^{-3}$. The substeps of the BUC scheme are integrated using explicit RK4. The singular value threshold for rank-adaptivity in the truncation step is set to $\vartheta = 10^{-5}$. In Figure~\ref{fig:3d3v_shear_rank_mass} we plot the rank and mass of the solution versus time. We observe that the rank remains relatively small for all time $t$ and the amount of mass lost is on the order of $10^{-3}$. 
In Figure~\ref{fig:3d3v_shear_vel} we plot slices in the $(x_2,x_3)$ plane of each component of the velocity with $x_1 \approx 0.49$ at time $t=5$ and $t=15$. 
In the 3d3v setting the vorticity has $3$ components 
\begin{equation} \label{eq:vort_3d}
\bm \omega = \left(
\partial_{x_2} u_3 - \partial_{x_3} u_2,
\partial_{x_3} u_1 - \partial_{x_1} u_3, 
\partial_{x_1} u_2 - \partial_{x_2} u_1
\right) 
\end{equation}
In Figure~\ref{fig:3d3v_shear_vor} we plot each component of the vorticity for the same slices.

\subsection{Explosion} 

\begin{figure}[!t]
\centering
\includegraphics[scale=0.26]{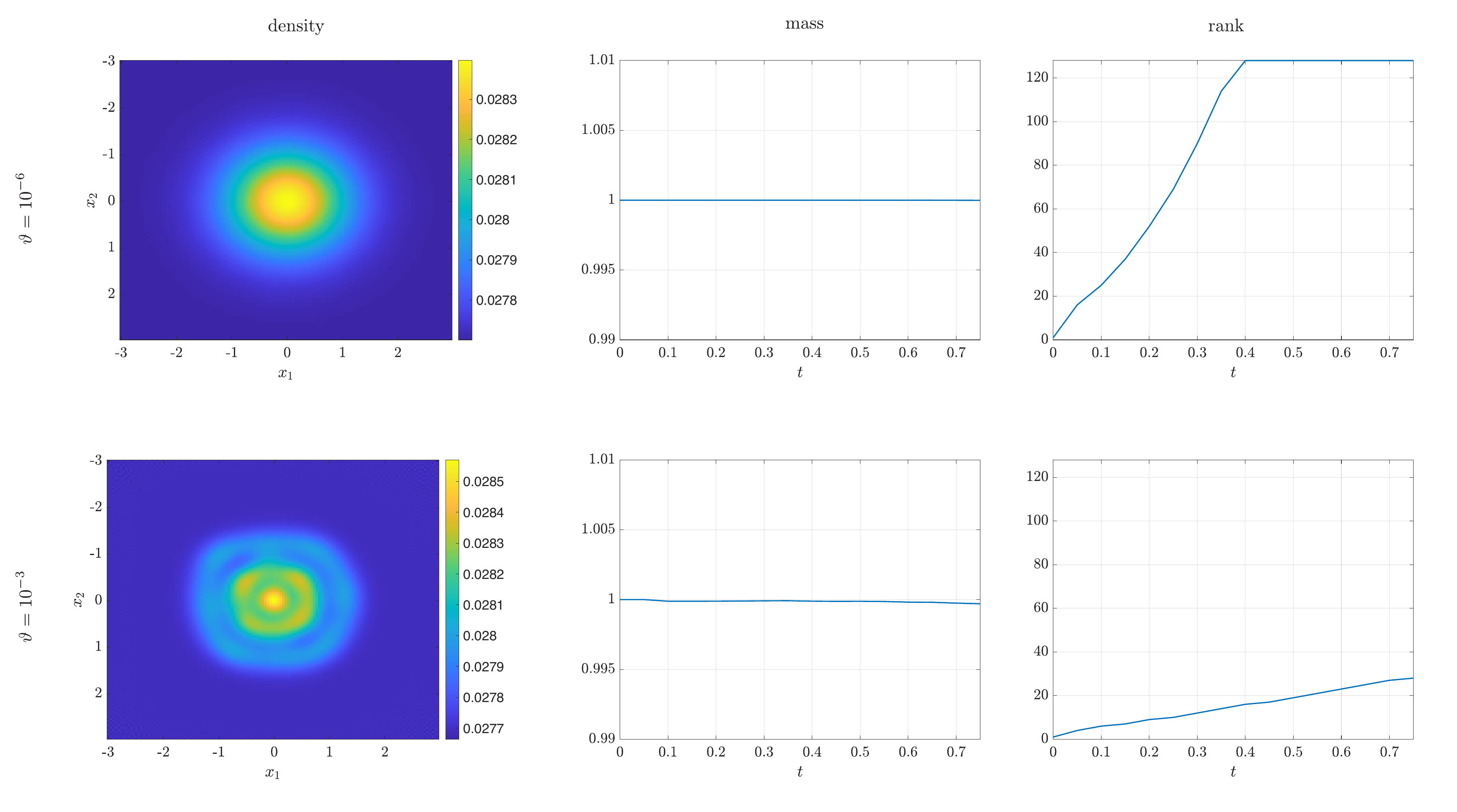}
\caption{Solution to the 2d2v explosion in the kinetic regime ($\epsilon=10$) computed with the BUC scheme. Results are shown for rank-adaptive threshold $\vartheta=10^{-6}$ and maximum rank $128$ (top row) and $\vartheta=10^{-3}$ (bottom row). The left column displays the density $\rho$ at time $t=0.75$, the middle column shows the mass versus time and the right column shows the rank versus time.} 
\label{fig:2d2v_kinetic}
\end{figure}

Next, we consider examples with large \text{density} initialized in a small region that then propagates outwards. To this end, we consider the rank-1 initial distribution function 
\begin{equation}
\begin{aligned} \label{eq:exp_IC} 
f(0,x,v) &= \frac{1}{M}X(0,x) V(0,v), \\ 
X(0,x) &= 1 + \alpha \prod_{i=1}^{d_x} \exp\left(\frac{-x_i^2}{2\sigma^2}\right) \\
V(0,v) &= \prod_{i=1}^{d_v}\exp\left(\frac{-v_i^2}{2}\right), 
\end{aligned}
\end{equation}
with $\alpha=0.25$, $\sigma=0.25$ and 
$$
%M = \int_{\Omega_x} 1 + \alpha \prod_{i=1}^{d_x} \exp\left(\frac{-x_i^2}{2\sigma^2}\right)dx \int_{\Omega_v} \prod_{i=1}^{d_v}\exp\left(\frac{-v_i^2}{2}\right)dv
M = \int_{\Omega_x} X(0,x) dx 
\int_{\Omega_v} V(0,v)dv. 
$$

\subsubsection{2d2v} 

First we set $d_x=d_v=2$ in \eqref{eq:exp_IC} and $\epsilon=10$ to demonstrate the interpolatory DLR method in the kinetic regime. We discretize the spatial domain $\Omega_x = [-3,3]^2$ using a tensor product grid with $128$ points per dimension and the velocity domain $\Omega_v = [-6,6]^2$ using a tensor product grid with $32$ points per dimension. We use the BUC scheme with time step-size $\Delta t = 10^{-3}$ and solved each substep with explicit RK4. We ran two simulations with rank-adaptive thresholds $\vartheta=10^{-3},10^{-6}$ and maximum rank $128$. In the left column of Figure~\ref{fig:2d2v_kinetic} we plot the density $\rho$ at time $t=0.75$, in the center column we plot the mass versus time and in the right column we plot the rank versus time for each simulation. For rank-adaptive threshold $\vartheta = 10^{-3}$ we observe that the solution rank does not exceed $30$ for $t \in [0,0.75]$ and the symmetry of the density is lost at $t=0.75$. Meanwhile, for rank-adaptive threshold $\vartheta=10^{-6}$ the solution rank reaches $128$ at $t \approx 0.4$ and the density at $t=0.75$ is symmetric. Our results indicate that, to obtain physically realistic solutions to the Boltzmann-BGK equation \eqref{eq:boltzmann-bgk} with the proposed interpolatory dynamical low-rank schemes, significantly larger ranks are required in the kinetic regime compared to the near-fluid regime. 

\begin{figure}[!t]
\centering
\includegraphics[scale=0.38]{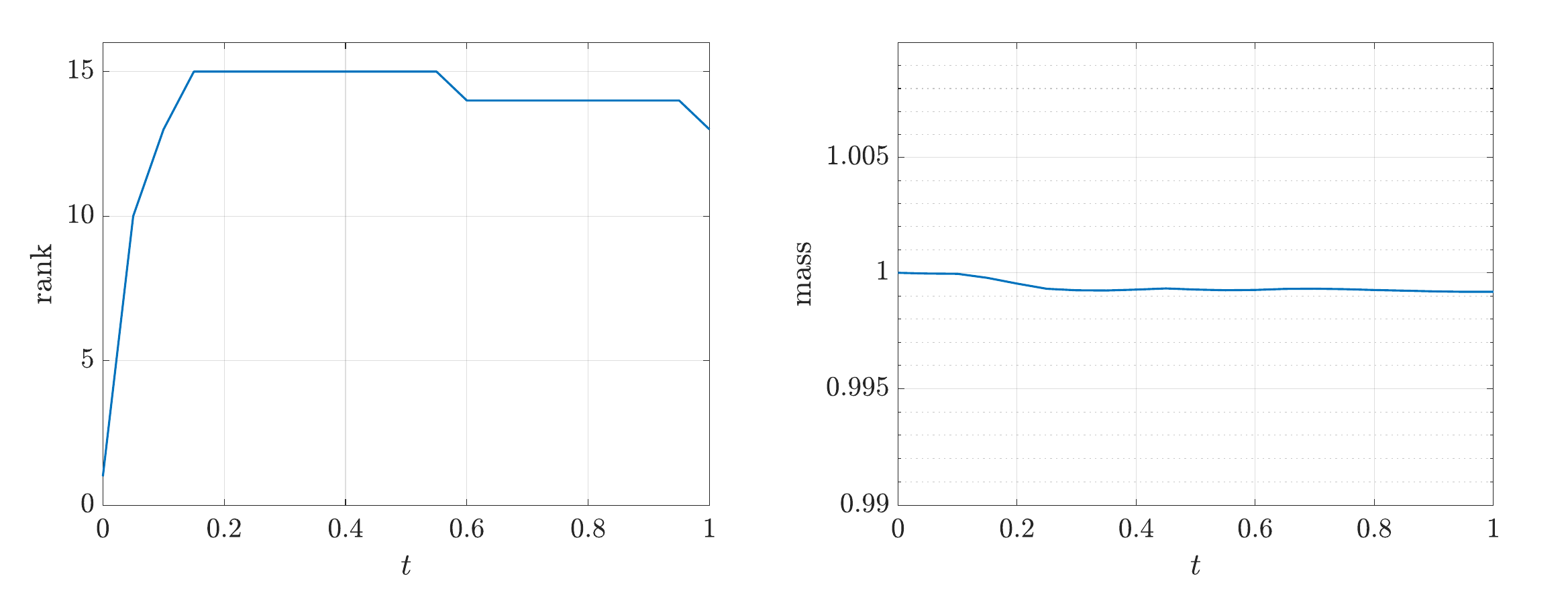}
\caption{Evolution of rank and mass over time for the BUC solution to the 3d3v Boltzmann-BGK equation, initialized with a rank-1 distribution given by \eqref{eq:exp_IC} and $\epsilon=10^{-3}$. } 
\label{fig:3d3v_ep1em3_rank_mass}
\end{figure}

\subsubsection{3d3v} 
\begin{figure}[!t]
\centering
\includegraphics[scale=0.38]{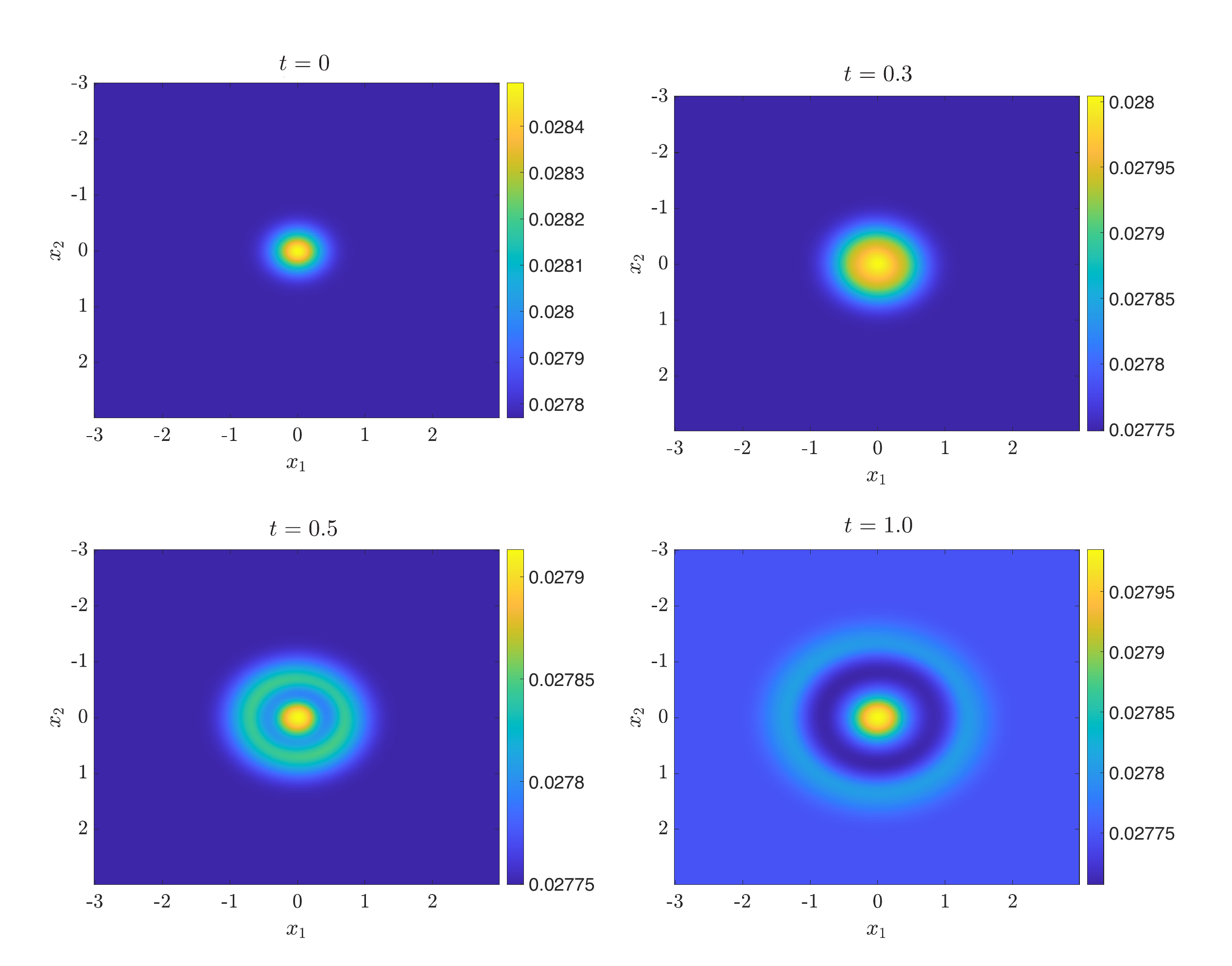}
\caption{$(x_1,x_2)$-marginal of the density for the BUC solution to the 3d3v Boltzmann-BGK equation initialized with rank-1 distribution given by \eqref{eq:exp_IC} and $\epsilon=10^{-3}$. Snapshots are shown at $t=0,0.3,0.5,1.0$. } 
\label{fig:3d3v_ep1em3_marg}
\end{figure}

Finally we set $d_x=d_v=3$ in \eqref{eq:exp_IC} and $\epsilon=10^{-3}$ to demonstrate the interpolatory DLR method in the fluid regime. We discretized the spatial domain $\Omega_x = [-3,3]^3$ using a tensor product grid consisting of 256 points per dimension and discretized the velocity domain $\Omega_v=[-6,6]^3$ using a tensor product grid consisting of 32 points per dimension. We used the BUC scheme with time step-size $\Delta t = 10^{-3}$ and solved each substep with explicit RK4 and rank-adaptive singular value threshold $\vartheta = 10^{-4}$. In Figure \ref{fig:3d3v_ep1em3_rank_mass} we plot the rank and mass of the low-rank BUC solution versus time. We observe that the rank remains smaller than $15$ for all time, which is significantly lower than the rank needed to represent the solution in the kinetic regime. In Figure \ref{fig:3d3v_ep1em3_marg} we plot time snapshots of the $(x_1,x_2)$ marginal functions at time $t=0.0,0.3,0.5,1.0$.

\section{Conclusions} \label{sec:conclusions}

We introduced two new solvers for the Boltzmann–BGK equation based on interpolatory dynamical low-rank (DLR) approximation. These methods extend the applicability of the DLR framework, enabling efficient time integration of the full {\color{black} six-dimensional (3d3v)} Boltzmann–BGK equation, which is unattainable with classical DLR techniques based on orthogonal projections. 
{\color{black} Numerical experiments in two-dimensional (1d1v) and four-dimensional (2d2v) dimensional settings, in both the fluid and kinetic regimes, demonstrate that the interpolatory DLR methods can achieve accuracy comparable to that of orthogonal DLR methods. We further demonstrated the interpolatory schemes on six-dimensional (3d3v) problems, where orthogonal DLR methods are prohibitively expensive.} 
Our results indicate that solutions in the kinetic regime may require significantly higher rank than those in the fluid regime. {\color{black} Nonetheless, the proposed interpolatory methods offer a promising direction for extending low-rank techniques to other nonlinear equations commonly arising in kinetic theory. }

\bibliographystyle{plain}
\bibliography{literature}

\end{document}